\documentclass[10pt,final]{siamltex}
\usepackage{times}
\usepackage{amsmath}
\usepackage{amsfonts}
\usepackage{amssymb}
\usepackage{graphicx}
\newcommand{\bfX}{\mathbf{X}}

\newcommand{\bfs}{\boldsymbol{\sigma}}
\newcommand{\refE}[1]{Equation~(\ref{#1})}

\newcommand{\refe}[1]{Eqn.~(\ref{#1})}

\newcommand{\refS}[1]{Section~\ref{#1}}
\newcommand{\refs}[1]{Sec.~\ref{#1}}

\title{
On hybrid simulation schemes for stochastic reaction dynamics }

\author{\small Sergey Plyasunov
\thanks{ Physics Department, UC Berkeley and
Physical Biosciences Division, E.O. Lawrence Berkeley National
Laboratory, 1 Cyclotron Road, Berkeley, CA 94720,e-mail:
teleserg@uclink.berkeley.edu } }

\begin{document}
\maketitle \centerline{\small current version:
 March 20,2005}
\begin{abstract}
The existing literature on stochastic simulation of chemical
reaction networks has a tendency to move as quickly as possible to
the abstract formulation of the stochastic dynamics in terms of
probabilities  based on the concept of the Chemical Master Equation
(CME), largely ignoring sample path  representation. In this
publication we discuss both theoretical basis and numerical approach
for the problems in this area using sample path methods as a crucial
part of the process.  Relying as it does on a representation of the
underlying stochastic processes as a weak solution of a system of
stochastic differential equations driven by Poisson random measures
this approach brings to bear a heretofore ignored but quite
effective problem solving methodology.
\par
 We  first present a simple and intuitive way of partitioning  species and reactions  of
the interaction network into different groups. We then discuss how
original stochastic dynamics with state dependent intensities of
transitions can be reformulated in terms of jump-diffusion
stochastic differential equations driven by both Wiener noise
sources and Poisson random measures. Finally, we show that this
approach facilitates the construction of hybrid simulation
techniques, an important step  in the creation of efficient
 techniques for modeling multi-scale stochastic  dynamics of the reaction networks.
 Numerical methods
related to sampling events from Poisson random measures are
demonstrated on simple intuitive examples. Error control analysis of
the  finite differences scheme is also presented.
\end{abstract}


\section{Introduction and Motivation\label{sec:intro}}
\par
 Many important problems in study of  complex systems
require consideration of  dense and nonlinear interaction between
different functional modules. Hybrid system is a framework first
proposed by the engineering community \cite{Sastry1999} to model
systems that exhibit both continuous and discrete state changes, and
has,  in recent years,  found increasingly wide applications in many
practical fields in connection to both complex man-made systems or
physical systems found in nature. However, most of the models
proposed in literature so far are deterministic, and are not
suitable to model inherent randomness.
\par
New data from high-throughput biology is enabling more ambitious,
complete and validated simulation biological models  and
 a lot of attention has recently been paid to the development of the
computational methods for the {\em in silico} investigation of complex biological processes on a
cellular level. To be valuable to biological and biomedical research computational methods must be
scaled to account for to large number of interacting component of different nature and  complex
pathways. As it was already noted by many researches stochastic effects play an important role in
the functionality of such systems \cite{ArkinMcAdamsRoss99} and approximation of molecular behavior
by continuous models fails to reproduce many interesting biological
phenomena\cite{RaoWolfArkin2002}. The deterministic approach  can be particularly misleading if the
molecular population of some critical reactant species becomes so small that microscopic
fluctuations can conspire with reaction channel feedback loops to produce macroscopic effects.
Recent work shows that this can happen with dramatic consequences in the genetic/enzymatic
reactions that go on inside a living cell \cite{RaoWolfArkin2002} or non-equilibrium combustion
reactions.
\par
 Correct methods for performing stochastic simulation of the discrete event
simulation are usually based on the Kinetic Monte Carlo Methods(KMC)
\cite{BortzKalosLebowitz75},Stochastic Simulation Algorithm(SSA)\cite{GillespieSSA,DTGillespieBook,
GibsonBruck2001}. KMC/SSA methods are exact i.e. they are designed to account for {\em every }
possible discrete reaction event and become notoriously inefficient for the scenarios when
stochastic effects are to be ignored and model can be abstracted as system of conventional ordinary
differential equations (ODEs) and ambitions goals which require the intensive use of Kinetic Monte
Carlo methods will also require immense computational resources. That became a main motivation for
the development of the hybrid simulation methods which will allow adequate levels of description
for  different parts of complex system in hand.
\par
Recently several attempts were made to bridge the gap between different levels of discretion and in
particular reduce excessive information content delivered by KMC/SSA methods on a temporal scale.
To overcome excessive level of resolution on the temporal scale i.e. to perform the {\em temporal
coarse-graining} $\tau$-leaping methods were developed
 \cite{GillespieTau,GillespiePerzoldI,Burage2004,RathinamGillespiePerzoldII}. Leaping  methods allow  to improve the
efficiency of  the standard KMC/SSA while maintaining acceptable losses in accuracy by
approximating the time-inhomogeneous Poisson process counting the reaction events via the Gaussian
process/diffusion approximation.
 The key idea here is to take a larger time steps and allow for more
  reactions to take place in that step, but under the proviso
that propensity functions do not  change too much in that interval. When it comes to practical
usage of the $\tau$-leaping method   a robust strategy is necessary  for  deciding when during the
coarse of the simulation when to step exactly using KMC/SSA and when to leap approximately with
diffusion methods \cite{TKurtz71,GillespieTau,EthierKurtz86}. When leaping, consideration of the
parameter values and error-reduction scheme is also necessary.
 \par
  This problem becomes of practical importance
 because it is often desirable to allow some chemical species to be
 treated  as continuous random variables and some to be treated
 discretely. This is particularly true for the case of the
 transcriptional regulation by transcription factors (TFs). In this
 case there may be as few as one DNA/TF binding site and messenger RNA (mRNA)
 abundance may be as small as tens of molecules
while  there may be  hundreds or even thousand
 molecules of  regulatory proteins per cell.
  In addition to that, biochemical networks, as probably every complex
system,  show remarkably interesting dynamics due to the presence of
   different time scales describing the evolution of discrete and
   continuous components. The technical difficulty with implementing the hybrid schemes is
   that standard KMC/Gillespie approach requires constant transition
   rates between reaction events  of the discrete species. This may
   not be the case if some of chemical species are evolving
   continuously in time or with different time scales.
  \par
 Several attempts have been made to illustrate the relevance and
 feasibility of hybrid algorithm, especially for the case of highly developed separation of time-scales.
 Rao and Arkin \cite{Rao2003}  consider the case of two groups of
 species evolving with two very
 different time scales and present  stochastic
 averaging procedure, constructing effective rates of the ``slow''
 reaction in the presence of the ``fast'' species, termed stochastic
 quasi steady state approximation (QSSA).
\par
In their work \cite{Kiehl2004}
 Kiehl et. al.  are developing an approach to integrate
 ODE/SDE schemes  with Kinetic Monte Carlo schemes.
 \cite{Kiehl2004} used synchronization algorithm  to produce the numerical
 scheme which bridges together asynchronous in nature KMC/SSA
 algorithms and  synchronous ODE integration schemes.
\par
 One should probably agree that following the   paths presented in the literature up to date
 will give solutions which  are  very close to the exact dynamics in the limit of a very small
 time step, largely different time scales, etc. but verity of different implementation
approaches presented in the literature to this point  indicates the fact that unifying solid
mathematical description has to be thought to justify the development of the hybrid simulation
schemes applicable for modeling stochastic phenomena in chemical reaction networks.
 Not many stochastic algorithms presented in the literature so far have gone
through the rigorous  error control analysis or even search for the fundamental justification of
the presented methods.
\par
\vskip 1cm In this publication we aim to provide not only the framework for the for the hybrid
simulation method but also we outline it's rigorous mathematical justification including error
analysis.
\par
Outline of this paper is following. In \refS{Notation::sec} we briefly outline some basic notation
and assumption used in this article. In the next  \refS{stoch::sec} our goal will be to give
precise mathematical formulation
 of the stochastic process corresponding to the hybrid representation of the reaction network.
We model the stochastic dynamics of the reaction network as a
dynamics of the  composite pair $(\mathbf{X},\boldsymbol{\sigma})$,
governed by the system of the jump-diffusion SDEs.
 This section  will be followed immediately by simple example and
view on error analysis. Article is ended with discussion and outline for the
future work.
\par
\section{Basic Notations and Assumptions\label{Notation::sec} }
\vskip 0.5cm
 Consider the reaction network which represent the
dynamics of $N$ well-stirred molecular species $\{
\mathsf{S}_i,i=1\ldots N \}$ which inter-react through $R$
($r=1,\ldots, R$) channels of chemical reactions. The current state
of the system is completely specified by the $N$-dimensional vector
$\mathbf{S}=(S_1,S_2,\ldots,S_N)$ consisting of non-negative integer
numbers, where each number $S_i$ defined to be the current total
number of molecules of type $\mathsf{S}_i$ in the system. Formally
reaction network can be expressed as a set of the transition rules:
\begin{subequations}\label{massaction}
\begin{gather}
r=1:~
\nu^{+}_{11}\mathsf{S}_1+\nu^{+}_{21}\mathsf{S}_2+\ldots+\nu_{N1}^{+}\mathsf{S}_N
\overset{k_1}{\rightarrow}
\nu^{-}_{11}\mathsf{S}_1+\nu^{-}_{21}\mathsf{S}_2+\ldots+\nu_{N1}^{-}\mathsf{S}_N,\\
r=2:~
\nu^{+}_{12}\mathsf{S}_1+\nu^{+}_{22}\mathsf{S}_2+\ldots+\nu_{N2}^{+}\mathsf{S}_N
\overset{k_2}{\rightarrow}
\nu^{-}_{12}\mathsf{S}_1+\nu^{-}_{22}\mathsf{S}_2+\ldots+\nu_{N2}^{-}\mathsf{S}_N,\\
\nonumber\dots\\
r=R:~
\nu^{+}_{1R}\mathsf{S}_1+\nu^{+}_{2R}\mathsf{S}_2+\ldots+\nu_{NR}^{+}\mathsf{S}_N
\overset{k_2}{\rightarrow}
\nu^{-}_{1R}\mathsf{S}_1+\nu^{-}_{2R}\mathsf{S}_2+\ldots+\nu_{NR}^{-}\mathsf{S}_N
\end{gather}
\end{subequations}
\par

 Transition rates  $a_r(\mathbf{S}): \mathbb{Z}^N\to \mathbb{R}_{+}$
 are the probabilities that reaction event corresponding to the
reaction $$r~:~\mathbf{S} \to\mathbf{S}  +{\boldsymbol \nu}_{r}$$
will take place in infinitesimal time interval $[t,t+dt)$ somewhere
inside the reaction volume.
Vectors ${\boldsymbol{\nu}_{r}}=\boldsymbol{\nu}_{r}^{-}-\boldsymbol{\nu}_{r}^{+}$ are the stoichiometric Al changes
of number of molecules of types $\mathsf{S}_i$ due to the  reaction $r$.
 As one can see, each reaction channel is specified by two quantities:
 transition rates $a_r$  and stoichiometric changes in molecular composition
of the system due to each individual event in
reaction channel $\boldsymbol{\nu}_r$.
\par
In \refe{massaction} $k_r$  are  transition  probabilities for
the reaction channel $r$ for any particular configurations of the
interacting molecules on the right hand side of (\ref{massaction}).
The functions $h_r(\mathbf{S}): \mathbb{Z}^{N}_{+}\to\mathbb{Z}$
are defined to be the number of distinct reaction combinations
of the reacting molecules available in the state $\mathbf{S}$.
The overall transition rate can be expressed as:
\begin{gather}\label{rate}
a_r(\mathbf{S})=k_rh_r(\mathbf{S}),
\end{gather}
 For the mass-action kinetics transition rates $h_r$ are found to
be a polynomial functions of variables $\mathbf{S}$.
For the simple reaction system:
$$\mathsf{S}_1\overset{k_1}{\to} \mathsf{S}_2,
~2\mathsf{S}_2\overset{k_2}{\to} \emptyset
 $$
statistical weights $h_i$ are:
$h_1(\mathbf{S})=S_1,~h_2(\mathbf{S})=\frac{1}{2}S_2(S_2-1)$.
 We refer  reader to the sources in
\cite{vanKampen92} and \cite{DTGillespieBook} for more details.
\section{Partition of species and reactions}
\vskip 0.5 cm
 Our interest in the next  step will be to consider
different reactions in the system in  a different manner. We
partition the set of molecular species into two groups: ``large
population group''
(with appropriate quantities labeled with "$X$" )
and discrete group (labeled with "$\sigma$"). Species from the first group  will be
represented by the set of $i=1,\ldots,N_X$ variables taking values
on $\mathbb{R}_{+}$:
$$
X_1,X_2,\ldots, X_{N_X},~~X_i\in~\mathbb{R}_{+}
$$
And present in large copy numbers: $X_i\gg 1,\forall i$
The reminder  of the $N_\sigma=N-N_X$ species represent the discrete
components of the network:
$$
\sigma_1,\sigma_2,\dots,\sigma_{N_\sigma},~~\sigma_i\in
\mathbb{Z}_{+}
$$
Partition of the species into the two separate groups correspond to
the partition of the stoichiometric vectors $\boldsymbol{\nu}_r$
into two sets :
$$
(\boldsymbol{\nu}^{X}_r,\boldsymbol{\nu}_r^{\sigma})
$$
which correspond to the changes of species in sets $\mathbf{X}$ and
$\mathbf{\sigma}$ due to the reaction event in the channel $r$.
\par
Based on this partition of species we can introduce three groups of
reactions:
\begin{romannum}
\item
 First group of reactions which we shall call $\mathcal{R}_1$, consist of
 reactions with large combinatorial weights, i.e. which satisfy the
 condition $h_r(\cdot)\gg 1$. We state that this condition is
 sufficient for the use of the diffusion approximation (see Appendix \refs{Appendix::diffusion})
\item
 Second group , $\mathcal{R}_2$ consists of reactions which do not
 satisfy the constraint $h_r\gg 1$ but change the species of the
 group $\mathbf{X}$, i.e. those reactions for which $\boldsymbol{\nu}^X_r\neq 0$
\item
Finally the third group, which we call $\mathcal{R}_3$ consists of
reactions for which $\boldsymbol{\nu}^{\sigma}_r\neq 0$
\end{romannum}
Note that, those subset of the reactions are not necessary
independent,i.e.
$$
\mathcal{R}_1\cap\mathcal{R}_3 \neq
\emptyset,~~\mathcal{R}_2\cap\mathcal{R}_3\neq \emptyset
$$
Without the loss of generality we shall assume that reactions from
the set $\mathcal{R}_2\cup\mathcal{R}_3$ (set of "jump" reactions)
are formed by reactions in \refe{massaction} with indices $1,\ldots,
R_d$, $R_d\leq R$. The rest,  $R-R_d$ are assumed to satisfy the condition of
having large combinatorial weight.
\section{Formulation of stochastic dynamics: ~
 jump-diffusion models\label{stoch::sec}}
\vskip 0.5cm \par
 Chemical Master equation
\cite{vanKampen92,handb_stoch_method} is usually employed to describe the stochastic dynamics of
reaction networks of the type (\ref{massaction}). It is  notoriously hard to solve analytically
even for the simple systems\cite{SamoilovCME} and, as we mentioned before, one usually has to
resort to so called {Kinetic Monte Carlo methods} \cite{GillespieSSA,BortzKalosLebowitz75} to
obtain numerical solution.
 In this article we shall abandon the description   based on the Chemical Master Equation(CME)
 and resort to the explicit treatment of the stochastic trajectories
 of the process $\{\mathbf{S}_t\}$ as a solution of the system of
stochastic differential equations (SDEs) providing
the recepie for the construction of  numerical schemes.\par
Poisson counter  process is a simple but important process and below we extend its
utility by combining with ideas from differential equations.
We briefly describe this  approach to pave our way to specific applications.
\par
To describe the structure of the stochastic jump processes let us
introduce in the rather traditional way  a stochastic basis
$(\Omega, \mathcal{F}, \mathbb{P})$ supplied with the filtration
 $\mathcal{F}_{t\geq 0 }$ large enough to include  $R$-dimensional
 unit  Poisson  processes. We recall that state variables $\mathbf{S}_t$ are defined on
  this  probability space. On this space transition rates \refe{rate} define the set of time non-homogeneous
 counting processes $N_r(t)$\cite{GihmanSkorohodSDE}:
\begin{subequations}
\begin{gather}
\mathbb{E}( N_r(t+\Delta t)-N_r(t)|\mathcal{F}_t )
=a_r(\mathbf{S}_{t_{-}})\Delta t+O(\Delta t^2),
\end{gather}
Notation $t_{-}$ is used to underline the existence of jumps in the process $\mathbf{S}_{t}$:
 $t_{-}=\lim_{\delta\to 0}(t-\delta)$. The above equation states that
differentials of  processes $N_r(\cdot)$, $dN_r=N_r(t+dt)-N_r(t)$ depend upon
local intensities $a_r(\mathbf{S}_{t_{-}})$. Notation $t_{-}$ is used to underline the existence of jumps in the process $\mathbf{S}_{t}$: $t_{-}=\lim_{\delta\to 0}(t-\delta)$.
 Counting process $N_r(\cdot)$ can be also
represented as a time changed  unit-rate Poisson processes $\Pi_r(\cdot),~r=1,\ldots, R$
\cite{EthierKurtz86}:
\begin{gather}
N_r(t)=\Pi_r\left(\int_0^t a_r(\mathbf{S}_{s_{-}})ds\right)
\end{gather}
\end{subequations}
In turn, counting processes $N_r(\cdot)$ represent the
stochastic dynamics of the Markov process $\mathbf{S}_t$  through the ``mass-balance'' equation:
\begin{gather}
\label{SDEPoisson}
\mathbf{S}_t=\mathbf{S}_0+\sum_{r=1}^R \boldsymbol{\nu}_rN_r(t),
\end{gather}
 Solution of \refe{SDEPoisson} is a constant on the interval where all  $N_r(\cdot)$ are constants and jumps by
 $\boldsymbol{\nu}_r$ whenever the reaction event takes place in the reaction channel
$r$ at the moment of time $t$:
\begin{gather}
\mathbf{S}_{t}=\mathbf{S}_{t_{-}}+\boldsymbol{\nu}_r
\end{gather}
Using the partitioning  which was introduced earlier,
 stochastic dynamics of the state vector $(\mathbf{X},\boldsymbol{\sigma})$
 can be presented as following system of  SDEs:
\begin{subequations}
\begin{gather}
\mathbf{X}_t=\mathbf{X}_0+\sum_{r=1}^R
\boldsymbol{\nu}^X_r N_r(t),\\
\boldsymbol{\sigma}_t=\boldsymbol{\sigma}_0+\sum_{r=1}^R \boldsymbol{\nu}^\sigma_rN_r(t),
\end{gather}
\end{subequations}
where integer vectors  $\boldsymbol{\nu}^X_r$ and
$\boldsymbol{\nu}^\sigma_r$ represent the change of the components
$\mathbf{X}$ and $\mathbf{\sigma}$ due to the reaction event in the
channel $r$. In general,  propensity functions $a_r(\cdot)$
 depend on both $\mathbf{X}$ and $\boldsymbol{\sigma}$.
\par
 Large combinatorial weights $h_r(\cdot)\gg 1$ of some reactions  allow us to follow the
{\em diffusion approximation} \cite{GillespieTau,EthierKurtz86},
i.e. by using the weak convergence methods we   express differentials of the time-inhomogeneous Poisson processes in \refe{SDEPoisson} corresponding the group $\mathcal{R}_1$ as a Gaussian processes with the drifts  $a_r(\cdot)=a_r(\mathbf{X}_{t_{-}},\boldsymbol{\sigma}_{t_{-}})$:
\begin{gather}
dN_r(\cdot)\rightarrow
a_r(\cdot)dt+\sqrt{a_r(\cdot)}dW_r(t),
\end{gather}
where $W_r(t),r\in \mathcal{R}_1$ are Brownian motions independent
of each other and $\mathbf{X}_{t=0}$.
This approximation of the stochastic dynamics  allows us to make
 a simplification in which we can consider the components $(\mathbf{X},\boldsymbol{\sigma})$
as defined on a hybrid space $\mathbb{R}^{N_X}_{+}\times\mathbb{Z}^{N_\sigma}_{+}$;
 $\mathbf{X}$ takes values on the continuum state space,
$\mathbb{R}^{N_X}_{+}$, while $\boldsymbol{\sigma}$ still takes values
on discrete state space$\mathbb{Z}^{N_\sigma}_{+}$.
 This brings us to the following formulation of the  stochastic dynamics of the pair
  $(\mathbf{X},\boldsymbol{\sigma})$
\begin{subequations}\label{jump_diffusion}
\begin{gather}\nonumber
d\mathbf{X}_t= \underbrace{\sum_{r\in \mathcal{R}_1}
\boldsymbol{\nu}^X_r
a_r(\mathbf{X}_{t_{-}},\boldsymbol{\sigma}_{t_{-}})dt+ \sum_{r\in
\mathcal{R}_1}  \boldsymbol{\nu}^\sigma_r
{a_r^{1/2}(\mathbf{X}_{t_{-}},\boldsymbol{\sigma}_{t_{-}})}dW_r(t)}_{diffusion}+\\
\label{jump_diffusionX} + \underbrace{\sum_{r\in
\mathcal{R}_2}dN_r(a_r(\mathbf{X}_{t_{-}},\boldsymbol{\sigma}_{t_{-}}))}_{jump},\\
\label{jump_diffusionSigma} d\boldsymbol{\sigma}_t=
\underbrace{\sum_{r\in\mathcal{R}_3}\boldsymbol{\nu}^\sigma_r
dN_r(a_r(\mathbf{X}_{t_{-}},\boldsymbol{\sigma}_{t_{-}}))}_{jump},
\end{gather}
\end{subequations}where we explicitly demonstrated the dependence of the
differential of the counting processes $dN_r(\cdot)$ on local
intensities. Note that according to the presence of the jump
component term in \refE{jump_diffusionX}
$\mathbf{X}_{t}\neq\mathbf{X}_{t_{-}}$ and $\mathbf{X}_t$ is a
jump-diffusion process in $\mathbb{R}^{N_X}_{+}$.
\par
Processes with state-dependent intensities of jump are notoriously hard to model. One complications
which makes the study of the SDEs of the type \refe{jump_diffusion} hard is the nature of the noise
source driving the dynamics. In contrast to very well studied Wiener or Poisson driven SDEs
\cite{KloedenPlaten92} is that intensity of the noise source is stochastic in nature. We shall take
a different representation of the noise source. Namely we use the technique of the Poisson random
measures. Our discussion of the  random measures will be brief and rather informal; for
mathematical background see monograph by Jacod and Shiryaev  \cite{JacodShiryaev} or Cox and Isham
\cite{CoxIshamBook}.
\par
 Poisson random measure $\mu(\omega,dt\times dz),\omega\in\Omega$ defines
a sequence of points and marks $\{(\tau_i,z_i)\}$ with the simple
interpretation that the mark $z_i$ arrives at time $\tau_i$;
~$\tau_i$ take values on $\mathbb{R}_{+}$. Marks $z_i$ take values
in a general space $E$ , which in the context of the problem
considered here can be taken as a subset of the Euclidian space,
$E\subset\mathbb{R}$.
 Deterministic intensity of the Poisson process $\mu(dt\times dz)$ is
$m(dz)~dt$ with $m(dz)$ being a  measure  on $E$. Arrival times
follow a Poisson process with deterministic intensity $dt$ and marks
a i.i.d distributed in the  interval $dz$.
\par
Important procedure which can be used to simulate the intensity dependent
process \refe{massaction} is based on introduction of the partition on the marked space.
 At every point of the state space $\mathbf{S}=(\mathbf{X},\boldsymbol{\sigma})$
we construct  the set of disjoint intervals based on transition rate
of "jump" reactions, i.e. reactions from the subset
$\mathcal{R}_2\cup\mathcal{R}_3$, i.e.  $r=1,\ldots,R_d$:
\begin{subequations}\label{intervals}
\begin{gather}
\Delta_{r}(\mathbf{S})=\left[\Lambda_{r-1}(\mathbf{S}),\Lambda_{r}(\mathbf{S})\right),\\
\Lambda_{r}=0,~\Lambda_{r}(\mathbf{S})=\sum_{r'=1}^{r}a_r(\mathbf{S})
\end{gather}
\end{subequations}
where index $r$ runs over reactions which either do not satisfy
condition $h_r \gg 1$ of for which vector
$\boldsymbol{\nu}^{\sigma}_r$ is non-zero i.e. they change discrete
components.  Thus in general length of the interval
$\Delta_r(\mathbf{S})$ is $a_r(\mathbf{S})$. Now define a set of
functions
$c_r:\mathbb{R}^{N_X}_{+}\times\mathbb{Z}_{+}^{N_\sigma}\times
\mathbb{R}\to \{0,1\}, ~r=1...R_d$:
\begin{gather}
c_r(\mathbf{X},\boldsymbol{\sigma},z)=1_{\Delta_r(\mathbf{S})}(z)=\left\{
\begin{array}{c}
1,~~{\rm ~if~} z\in\Delta_r(\mathbf{X},\boldsymbol{\sigma}),\\
 0,~~{\rm otherwise}
\end{array}
\right.
\end{gather}
Suppose now that the space $(\Omega, \mathcal{F},\mathcal{F}_t,P)$
supports a Poisson random measure $\mu(\omega,ds\times
dz)=\mu(dt\times dz),~\omega\in\Omega$ in $\mathbb{R}_{+}\times
E,~E\subset\mathbb{R}$ with a compensator $dt~dz$ and a set of
Wiener processes $W_r(\cdot),~r\in \mathcal{R}_1$. Then time
non-homogenous counting process $N_r(t)=\int_0^t
dN_r(a_r(\mathbf{X}_{s_{-}},\boldsymbol{\sigma}_{s_{-}}))$ can be
represented as \cite{GihmanSkorohodSDE,PProtter}:
\begin{gather}
N_r=\int_0^t\int_{\mathbb{R}}
c_r(\mathbf{X}_{s_{-}},\boldsymbol{\sigma}_{s_{-}},z) \mu(ds\times
dz)
\end{gather}
 and  $(\mathbb{R}_{+}^{N_c}\times \mathbb{Z}_{+}^{N_d})$ valued
process $(\mathbf{X}_t,\boldsymbol{\sigma})$ can be viewed as a
solution of coupled Ito-Skorohod SDEs driven by the Poisson random
measures \cite{PProtter,GihmanSkorohodSDE}:
\begin{subequations}
\label{random_measure}
\label{jump_diffusion_poisson}
\begin{gather}\nonumber
d\mathbf{X}_t=\sum_{r\in \mathcal{R}_1} \boldsymbol{\nu}^X_r
a_r(\mathbf{X}_{t_{-}},\boldsymbol{\sigma}_{t_{-}})dt+ \sum_{r\in
\mathcal{R}_1}  \boldsymbol{\nu}^\sigma_r
{a_r^{1/2}(\mathbf{X}_{t_{-}},\boldsymbol{\sigma}_{t_{-}})}dW_r(t)+\\+
\sum_{r\in \mathcal{R}_2}\int_{\mathbb{R}}\boldsymbol{\nu}^\sigma_r
c_r(\mathbf{X}_{t_{-}},\boldsymbol{\sigma}_{t_{-}},z) \mu(dt\times dz),\\
d\boldsymbol{\sigma}_t=\sum_{r\in\mathcal{R}_3}
\int_{\mathbb{R}}\boldsymbol{\nu}^\sigma_r
c_r(\mathbf{X}_{t_{-}},\boldsymbol{\sigma}_{t_{-}},z)
 \mu(dt\times dz)
\end{gather}
\end{subequations}
 Function $c_r$ in \refe{random_measure} transforms arrived marks into magnitude of
 the
jump of $\mathbf{S}$ due to appropriate reaction channel.
\par
For the practical insight on the dynamics described above, we will introduce {\em reference Poisson process}.
\par
 Since the length of each interval $\Delta_r(\mathbf{S})$ is $a_r(\mathbf{S})$
and bounded function for every $r=1\ldots R$ then it is follows that
the length of the interval $\Delta(\mathbf{S})$ (we denote it
$|\Delta(\mathbf{S})|$) is bounded  function of
$\mathbf{S}=(\mathbf{X},\boldsymbol{\sigma})$. Therefore, it has a
maximum at some point $S^*=(\mathbf{X}^*,\boldsymbol{\sigma}^*)$:
\begin{gather}
|\Delta(S)|\leq |\Delta(S^*)|
\end{gather}
Let $\Lambda_{max}=|\Delta(S^*)|$ denote the maximum length of the
joint intervals in the marked space. Now  consider  a standard
Poisson process with intensity $\Lambda_{max}$. Let
$\tau^{\Lambda}_n,n=1,2,\ldots$ denote the jump times of this
process (sampled, perhaps, as the partial sums  of the exponentially
distributed random variables with mean $\Lambda_{max}^{-1}$). Let
$\Delta^*=\Delta(\mathbf{S}^*)$ be the marked space and $\lbrace
z_n\rbrace$ be the sequence of i.i.d. distributed random variables
with uniform distribution on $\Delta^*$, independent on $N_t$. In
this case we can represent a Poisson random measure with intensity
$dz~dt$  is  related to the marked point process
$(\tau^{\Lambda}_n,z_n)_{n\geq 0}$, constituting the sample path of
the process. Namely for each $A\in \Delta^*$:
\begin{subequations}
\begin{gather}\label{Poisson_master}
N_{\Lambda}((0,t],A)=\sum_{n\geq 1}1_{ \tau^{\Lambda}_n\leq t } 1_{z_n \in A}, \\
\mathbb{E}(N_{\Lambda}((0,t],A))=\Lambda_{max}t\mathbb{P}(z_n\in A)=
t\mu(A)
\end{gather}
\end{subequations}
This representation (sometimes referred to as thinning of the
Poisson measure \cite{CoxIshamBook}) is very useful in practical
numerical applications. In the next section we will use the
reference process $N_{\Lambda}$ to construct the   discretization
scheme taking $A=\Delta_r(\cdot)$ (\ref{intervals})
 for various $r\in\mathcal{R}_{2,3}$.
\section{Implementation of Numerical Scheme}
\par In this section we outline basic principles for the
construction of the numerical scheme for the solution of the system
given by \refe{jump_diffusion_poisson}. To begin with, we should
introduce the appropriate discretization of the interval $[0,T]$ on
which dynamics of the system has to be investigated. Let us denote
by $[0,T]_{h,M}=\lbrace 0=t_0<t_1<t_2<\ldots<t_{M+1}=T \rbrace$,
perhaps the usual equidistant, discretization of the time interval ,
i.e. $h=T/M$. \par
 Suppose now that $\lbrace
0,\tau_1^{\Lambda},\tau_2^{\Lambda},\dots\rbrace,~~\tau_i^{\Lambda}\in
[0,T]$ are jump times of the reference Poisson process
\ref{Poisson_master} ; then we consider new time discretization
which is a merger of the points from  $[0,T]_{h,N}$ and jumps of the
reference Poisson process \ref{Poisson_master}.
\par
It is important to note that jump-times of the Poisson measure
\refe{Poisson_master} can be modeled without discretization error
and with little prior information; one has only make sure that
transition rates from the groups $\mathcal{R}_{2,3}$ are bounded,
i.e. are less then chosen intensity of the reference Poisson
measure. Armed with the time discretization scheme assembled from
equidistant  discretization and jump-times of the reference process
we introduce the  following simple jump-diffusion numerical scheme:
\begin{subequations}
\begin{gather}\label{Euler_diff_X}\nonumber
\hat{X}_{j,t_i}^{-}=\hat{X}_{j,t_{i-1}}+\sum_{r\in\mathcal{R}_1}\nu^{X}_{jr}
a_r(\hat{\mathbf{X}}_{t_{i-1}},\hat{\boldsymbol{\sigma}}_{t_{i-1}})(t_i-t_{i-1})+\\+
\sum_{r\in\mathcal{R}_1}
\nu^{X}_{jr}a^{1/2}_r(\hat{\mathbf{X}}_{t_{i-1}},\hat{\boldsymbol{\sigma}}_{t_{i-1}})(W_r(t_i)-W_r(t_{i-1})),\\
\label{Euler_jumps_sigma}
\hat{X}_{j,t_i}=\hat{X}^{-}_{j,t_{i}}+\sum_{r\in\mathcal{R}_2}\nu^{X}_{jr}\int_{\Delta^*}
c_r(\hat{\mathbf{X}}^{-}_{t_i},\hat{\boldsymbol{\sigma}}_{t_{i-1}},z)\mu(dt\times dz),\\
\hat{\sigma}_{j,t_i}=\hat{\sigma}_{j,t_{i-1}}+
\sum_{r\in\mathcal{R}_3}\nu^{\sigma}_{jr}\int_{\Delta^{*}}
c_r(\hat{\mathbf{X}}^{-}_{t_i},\hat{\boldsymbol{\sigma}}_{t_{i-1}},z)~\mu(dt\times
dz)
\end{gather}
\end{subequations}
which recursively determines the values of the  discretized
processes $\hat{\mathbf{X}}=(\hat{X}_1,\ldots,\hat{X}_{N_c})$ and
$\hat{\boldsymbol{\sigma}}=(\hat{\sigma}_1,\ldots,\hat{\sigma}_{N_d})$
at points $t_i$, starting from values $X_0,\sigma_0$. Increments of
the Winer processes $W_r(\cdot)$ are zero-mean Gaussian random
variables: $W_r(t_i)-W_r(t_{i-1})\propto
\mathcal{N}(0,t_i-t_{i-1})$.
 According to the (\ref{Euler_diff_X}),
 between the jumps dynamics of the component $\mathbf{X}$ is purely diffusive.
  Integral in \refe{Euler_jumps_sigma} is
evaluated at the single point $z_i$ which belongs to the uniform
distribution $\mathcal{U}([0,\Lambda_{max}])$, i.e. jump in channel
$r$ is  accepted iff generated random variable
$z_{t_i}\in\mathcal{U}([0,1])$, generated on the step $t_i$
satisfies the condition:
\begin{gather}
z_{t_i}\Lambda_{max}\in
\Delta_r((\hat{X}_{j,t_i}^{-},\hat{\boldsymbol{\sigma}}_{t_{i-1}})),
\end{gather}
\par
Depending on the relative ratio between characteristic time-scales
of diffusion and discrete species different numerical schemes must
be considered to achieve substantial  speed-up of the simulation.
Here we consider only the case of fast diffusion modes and slow
discrete modes. In this case every interval
$[\tau^{\Lambda}_{n},\tau^{\Lambda}_{n+1})$ should be partitioned
onto smaller intervals of the discrete grid on which numerical
scheme (\ref{Euler_diff_X}) for diffusion approximation is used.
Opposite case of the fast switching in discrete component $\bfs$
will be considered elsewhere.
\section{Example}
We now wish do demonstrate application of the above numerical scheme.
In this section we consider the simple example, which nevertheless, is interesting enough
for the purposes of the demonstration of the technique discussed in this paper.
Consider the following model of the reaction network:
\begin{subequations}
\begin{gather}
\mathsf{S}_1\overset{k_1}{\rightarrow }\mathsf{S}_2,\\
\mathsf{S}_2\overset{k_2}{\rightarrow }\mathsf{S}_1,\\
\mathsf{S}_1\overset{k_3}{\rightarrow }\mathsf{S}_1+n\mathsf{S}_3,\\
\mathsf{S}_2+\mathsf{S}_3\overset{k_4}{\rightarrow}\mathsf{S}_2+\mathsf{S}_4,\\
\mathsf{S}_4\overset{k_5}{\rightarrow}\emptyset
\end{gather}
\end{subequations}
with continuous and discrete components $\boldsymbol{\sigma}=(S_1,S_2)=(\sigma_1,\sigma_2)$ and
$\mathbf{X}=(S_3,S_4)=(X_1,X_2)$ respectively. The functions $a_r$ and vectors $\boldsymbol{\nu}_r$
for this system are:
\begin{subequations}
\begin{gather}
a_1(\mathbf{X},\boldsymbol{\sigma})=k_1\sigma_1,~ \boldsymbol{\nu}_1^T=(-1,1,0,0),\\
a_2(\mathbf{X},\boldsymbol{\sigma})=k_2\sigma_2,~ \boldsymbol{\nu}_1^T=(1,-1,0,0),\\
a_3(\mathbf{X},\boldsymbol{\sigma})=k_3\sigma_1,~ \boldsymbol{\nu}_1^T=(0,0,n,0),\\
a_4(\mathbf{X},\boldsymbol{\sigma})=k_4\sigma_2X_1,~ \boldsymbol{\nu}_1^T=(0,0,-1,1),\\
a_5(\mathbf{X},\boldsymbol{\sigma})=k_5X_2,~ \boldsymbol{\nu}_1^T=(0,0,0,-1)
\end{gather}
\end{subequations}
We have used the following set of kinetic parameters:$k_1=0.50~~,k_2=0.50~,k_3=1.00~~,k_4=0.10~~$,
$k_5=0.01~~$, $X_1(0)=1000,X_2(0)=200$ and $n=5$\footnote{ One can consider this as a simple model
of transcriptional regulation, where presence/absence of the transcription factor $(S_1\in
\{0,1\})$ leads to the bursts of the transcription of protein $S_3$ with $n$ proteins per
transcription event.}.
\par
 Analysis of the propensity functions ($h_{4,5}\gg 1$) shows the following
partition of the original reaction set $\mathcal{R}=\{1,\ldots,5\}$:
\begin{gather}
\mathcal{R}_1=\lbrace 4,5\rbrace,\\
\mathcal{R}_2=\lbrace 1,2\rbrace,\\
\mathcal{R}_3=\lbrace 3\rbrace
\end{gather}
i.e. reactions $\lbrace 4,5\rbrace$ correspond to the diffusion modes.
 We also constraint ourselves in this example considering  the case:
\begin{gather}
a_1+a_2+a_3=\Lambda_{max} < a_{4,5}
\end{gather}
 Propagation of the diffusion mode was performed via Euler-Maruyama
 scheme \cite{Burage2004},\cite{KloedenPlaten92} with set of different
time-steps: $h=0.1,\ldots,2.0$. We performed $10^4$ Monte Carlo runs with both standard KMC/SSA and
hybrid scheme outlined above to obtain the distributions of the number of molecular species
$(S_3,S_4)=(X_1,X_2)$ at time-slice $T=2\times 10^3$.
 \par
  Figures(\ref{Fig::histogramX1},\ref{Fig::histogramX2}) show
  that distributions of the species $S_4,S_5$, $P(S_3,T)=P(X_1,T), P(S_4,T)=P(X_2,T) $
   obtained with the hybrid scheme is almost indistinguishable from the results of the
 obtained by following complete KMC/SSA-based approach. We also report
 an observed speed-up in terms  of the ratio of the CPU times $\frac{T_{KMC}}{ T_{hybrid}}$,
 Fig.(\ref{Fig::speed_up}) as a function of time step $h$ used in numerical  integration of diffusion modes.
\vskip 0.5cm
\section{Discussion}
\par
Advances in field of molecular biology outlined new  directions of
research in stochastic chemical kinetics. We now need to develop the
software tools and mathematical approaches to integrate models from
micro-scales to macro-scales in a seamless fashion. Such multi-scale
models are essential if we are to produce quantitative, predictive
models of complex biological behaviors. An important step on this
path is implementation of  hybrid simulation methods, capable to
account for heterogeneity of properties of interacting components.
\par
We have outlined rigorous framework for development of hybrid
simulation schemes based on the path sample representation of
reaction dynamics rather then on CME based approach. For the best of
our knowledge this type of analysis has not being  performed before.
\par It was not the primal goal of this publication to justify the partitioning of the
species/reactions into the different groups. One should probably
have some prior knowledge when partitioning the species into
different types. Choice of the reactions ("diffusive" or "jump") is
based on functional law of large numbers and  error generated by
this step can probably be controlled rigorously with some insights
on this subject outlined in \refS{Appendix::diffusion}.
\par
Numerical solution for the propagation of the diffusion modes were outlined via simple constant
step explicit Euler scheme. Typical explicit or implicit scheme \cite{RathinamGillespiePerzoldII}
based on stochastic Taylor approximation \cite{KloedenPlaten92},\cite{Milstein1995} might run into
problems , since they do not conserve the non-negativity of the numerical solution, i.e. i.e. it
 cannot guarantee $\hat{X}_{j,t_i}\in \mathbb{R}_{+}$ almost surely.
 We point out here that use of balanced-implicit stochastic schemes (BIM),
 \cite{MillsteinPLatenSchurz98} look promising when it comes to the construction of the
 numerical solution with the property of ``almost sure'' positivity.
BIM scheme known to have the same order of convergence as the
Euler-Maruyama stochastic scheme , namely , the error $O(h^{1/2})$
for approximations of the individual trajectories and $O(h)$ for
moments.
 \par
Adaptive time-step control systems also look promising as a
direction for the research aimed for achieving the simulation speed
up\cite{Burage2004}.
\par

\section{ACKNOWLEDGMENTS}
Author  would like to acknowledge valuable discussions and useful comments from Prof. C. Myers,
Prof. D. Steinsaltz, K. Ericksson and Prof. A.P. Arkin.
\appendix{
\section{Diffusion limit of the Poisson process with state-dependent intensity
\label{Appendix::diffusion}}
\par
 Introducing the set of independent Poisson processes with rate
$k_r$ differential of the counting process $dN_r$  at some state
$\mathbf{S}_t$ can be expressed as:
\begin{gather}
dN_r=\sum_{i=1}^{n=h_r(\mathbf{S})}d\Pi_{k_r}^{(i)}=
k_rh_r(\mathbf{S})dt+dM_r,\\
dM_r=\sum_{i=1}^{n=h_r(\mathbf{S})}d\tilde{\Pi}_{k_r}^{(i)},
~d\tilde{\Pi}_{k_r}^{(i)}=d\Pi_{k_r}^{(i)}-k_rdt,~\mathbb{E}(d\tilde{\Pi}^{(i)}_{k_r}|\mathcal{F}_t)=0
\end{gather}
Using the functional analog of the law of the large numbers one can
conclude:
\begin{gather}
dM_r=\sqrt{k_rh_r}dW_r+d\epsilon(t),~~
\end{gather}
where  $\mathbb{E}(\|d\epsilon_r(t)\|~\|\mathcal{F}_t)\propto
h_r^{-1}(\mathbf{S})\ll 1$ as $h_r\gg 1$ (see
\cite{EthierKurtz86,Anantahram95}). Note that this condition is
different from previously reported in the literature
\cite{GillespiePerzoldI}.

\section{Convergence analysis of the Hybrid Approximation}
\par
Below we demonstrate convergence properties of the numerical scheme
(\ref{Euler_diff_X},\ref{Euler_jumps_sigma}) if propensity functions
$a_r(\cdot)$ satisfy regular Lipschitz continuity  condition:
\begin{gather}\label{Lipshitz_cond}\nonumber
|a_r(\bfX,\bfs)-a_r(\bfX',\bfs')|^2+|a_r^{1/2}(\bfX,\bfs)-a_r^{1/2}(\bfX',\bfs')|^2\leq
\\\leq K(|\bfs-\bfs'|^2+|\bfX-\bfX'|^2),~~\forall r\in\mathcal{R},
\end{gather}
where  $|\bfX|$,$|\bfs|$ are  the Euclidian norms of the vectors
$\bfX$,$\bfs$ respectively  and $K $ is some constant.
\par
Since in  scheme \ref{jump_diffusionSigma} generation of the jump
times of the counting process $N_{\Lambda}(t)$ is not linked to the
dynamics of the state space pair $(\mathbf{X},\boldsymbol{\sigma})$
we can consider first the error of the approximation of the scheme
(\ref{jump_diffusionSigma}) for a given sequence of the jump times
of the Poisson process (\ref{Poisson_master}) $N_{\Lambda}$ :
$$
0<\tau^{\Lambda}_1<\tau^{\Lambda}_{2}<\ldots<\tau^{\Lambda}_{N}<T,~~
N_{\Lambda}(T)=N,
$$
with random number $N$ following the Poisson distribution:
\begin{gather}
P(N_{\Lambda}(T)=N)=\frac{(\Lambda_{max}T)^N}{N!}\exp(-\Lambda_{max}T)
\end{gather}
\par
On the interval $[\tau_{k}^{\Lambda},\tau_{k+1}^{\Lambda})$ we
represent the numerical solution corresponding to the diffusion mode
as:
\begin{gather}\label{Appendix::scheme}
\hat{\bfX}_t=\hat{\bfX}_{\tau_{k}^{\Lambda}} +
\int_{\tau_{k}^{\Lambda}}^t \sum_{r\in\mathcal{R}_1}\nu_r
{a}_r(\hat{\bfX}_s,\hat{\bfs}_{\tau_{k}^{\Lambda}})ds+
\sum_{r\in\mathcal{R}_1}\nu_r{a}^{1/2}_r(s,\hat{\bfX}_s,\hat{\bfs}_{\tau_{k}^{\Lambda}})dW_r,\\
\hat{\bfX}_t=\hat{\bfX}_{t_k},~t\in [t_k,t_{k+1})\subset
[\tau_{k}^{\Lambda},\tau_{k+1}^{\Lambda}),~
~\hat{\bfX}_{t=0}=\bfX_0,~\hat{\bfs}_{t=0}=\bfs_0.
\end{gather}
 On each interval between two jump times of the reference process
$[\tau_{n}^{\Lambda},\tau_{n+1}^{\Lambda})$
 we perform our analysis in a rather traditional way, starting from
estimating the expectation of the norms $ |{\bfX}_t-\hat{\bfX}_t|^2,
 |{\bfs}_t-\hat{\bfs}_t|^2
$
conditional on filtration generated up to  jump time
$\tau^{\Lambda}_n$: $\mathcal{F}_{\tau^{\Lambda}_n}$:
\begin{subequations}
\begin{gather}\label{error_norm:SDE}
\epsilon_X([\tau_n,\tau_{n+1}))=
 \sup_{t\in[\tau_{n}^{\Lambda},\tau_{n+1}^{\Lambda})}\mathbb{E}(
|{\bfX}_{t}-\hat{\bfX}_{t}|^2 ~|\mathcal{F}_{\tau^{\Lambda}_{n}}),\\
\epsilon_\sigma([\tau_n,\tau_{n+1}))=
 \sup_{t\in[\tau_{n}^{\Lambda},\tau_{n+1}^{\Lambda})}\mathbb{E}(
|{\bfs}_{t}-\hat{\bfs}_{t}|^2 ~|\mathcal{F}_{\tau^{\Lambda}_{n}})
\end{gather}
\end{subequations}
Then, evidently, total error on the interval $[0,T]$ for each of the
components, $\bfX$ and $\bfs$ can be found as:
\begin{gather}
\epsilon_X([0,T))=\max_n\epsilon_X([\tau_n,\tau_{n+1})),\\
\epsilon_\sigma([0,T))=\max_n\epsilon_\sigma([\tau_n,\tau_{n+1}))
\end{gather}
Let us turn to the estimation of (\ref{error_norm:SDE}) on each
interval $[\tau_{n}^{\Lambda},\tau_{n+1}^{\Lambda})$ based on the
formulation (\ref{Appendix::scheme}). For the component $\bfX$ for
every $t\in [\tau_{n}^{\Lambda},\tau_{n+1}^{\Lambda})$ one has:
\begin{gather}\nonumber
\epsilon_X([0,t])=\sup_{s\in[0,t]}\mathbb{E}( [
\hat{\bfX}_{\tau_n}-\bfX_{\tau_n} +\int_0^s
\sum_{r\in\mathcal{R}_1}\nu_r
(a_r(\bfX_u\bfs_u)-a_r(\hat{\bfX}_u,\hat{\bfs}_u))du+\\+
\sum_{r\in\mathcal{R}_1} \nu_r
(a_r^{1/2}(\bfX_u\bfs_u)-a_r^{1/2}(\hat{\bfX}_u,\hat{\bfs}_u))dW_r(u)
 ]^2|\mathcal{F}_{\tau_n^{\Lambda}} )\leq \\ \leq
\label{Appendix:rhs:error}
 |\hat{\bfX}_{\tau_n}-\bfX_{\tau_n}|^2+ \sup_{s\in
[0,t]}\mathbb{E}(|A_1(s)|^2|\mathcal{F}_{\tau_n^{\Lambda}})
+\sup_{s\in
[0,t]}\mathbb{E}(|A_2(s)|^2|\mathcal{F}_{\tau_n^\Lambda})
\end{gather}
where:
$$
A_1(s)=\int_{\tau_n}^s\sum_{r\in\mathcal{R}_1}\nu_r
(a_r(\bfX_u,\bfs_{\tau_n})-a_r(\hat{\bfX}_u,\hat{\bfs}_{\tau_n}))du,
$$
and
$$
A_2(s)=\int_{\tau_n}^s \sum_{r\in\mathcal{R}_1} \nu_r
(a_r^{1/2}(\bfX_u,\bfs_{\tau_n})-a_r^{1/2}(\hat{\bfX}_u,\hat{\bfs}_{\tau_n}))dW_r(u)
$$
By Couchy's inequality  bound for the  second term in
\ref{Appendix:rhs:error} can be presented as:
\begin{gather}
\sup_{s\in
[0,t]}\mathbb{E}(|A_1(s)|^2|\mathcal{F}_{\tau_n^{\Lambda}})\leq
(t-\tau_n^{\Lambda})\sup_{s\in [0,t]}\mathbb{E}(
\int_{\tau_{n}^{\Lambda}}^s (\sum_{r\in\mathcal{R}_1} \nu_r
(a_r(\bfX_u,\bfs_{\tau_n})-a_r(\hat{\bfX}_u,\hat{\bfs}_{\tau_n})))^2du
 |\mathcal{F}_{\tau_n^{\Lambda}}),
\end{gather}
At the same time, by Doob's inequality  for martingales, third term
can be estimated as:
\begin{gather}
\sup_{s\in
[0,t]}\mathbb{E}(|A_2(s)|^2|\mathcal{F}_{\tau_n^{\Lambda}})\leq 4
\mathbb{E}(|A_2(t)|^2|\mathcal{F}_{\tau_n^{\Lambda}}) = 4
\mathbb{E}\int_{\tau_n^{\Lambda}}^t\sum_{r\in\mathcal{R}_1} \nu_r
(a_r^{1/2}(\bfX_u,\bfs_{\tau_n})-a_r^{1/2}(\hat{\bfX}_u,\hat{\bfs}_{\tau_n}))^2du
\end{gather}
which together with Lipschitz condition (\ref{Lipshitz_cond}) and
Gronwall lemma gives the bound for the error of the component $\bfX$
on interval $[\tau_n^{\Lambda},t]$,
$\epsilon_X([\tau_{n}^{\Lambda},t])$:
\begin{subequations}
\begin{gather}
\epsilon_X([\tau_{n}^{\Lambda},t])\leq K\int_{\tau_n^{\Lambda}}^t
\epsilon_{X}([\tau_n^{\Lambda},s])ds +
KT\mathbb{E}(|\hat{\bfX}_{\tau_n^{\Lambda}}-\bfX_{\tau_n^{\Lambda}}|^2
+KT\mathbb{E}(|\bfs_{\tau_n^{\Lambda}}-\hat{\bfs}_{\tau_n^{\Lambda}}|^2),\\
\label{Appendix::Granwall} \epsilon_X([\tau_{n}^{\Lambda},t])\leq
KT[\mathbb{E}(|\hat{\bfX}_{\tau_n^{\Lambda}}-\bfX_{\tau_n^{\Lambda}}|^2+
KT\mathbb{E}(|\bfs_{\tau_n^{\Lambda}}-\hat{\bfs}_{\tau_n^{\Lambda}}|^2)]\exp(K(t-\tau_n^{\Lambda}))
\end{gather}
\end{subequations}
One can see that additional error due to the possible jump at time
$\tau_{n+1}^{\Lambda}$ terms is also bounded:
\begin{gather}\label{Appendix::jump}
\int_{\mathbb{R}}
|c_r(\bfX_{\tau_{n+1-}^{\Lambda}},\bfs_{\tau_{n+1-}^{\Lambda}},z)-
c_r(\hat{\bfX}_{\tau_{n+1-}^{\Lambda}},\hat{\bfs}_{\tau_{n+1-}^{\Lambda}},z)|^2dz
\leq
K(|\hat{\bfX}_{{\tau_{n+1-}^{\Lambda}}}-\bfX_{\tau_{n+1-}^{\Lambda}}|^2+
|\hat{\bfs}_{\tau_{n+1-}^{\Lambda}}-\bfs_{\tau_{n+1-}^{\Lambda}}|^2),
\forall r\in \mathcal{R}_{2,3},
\end{gather}
\par
Assuming that one is using scheme with $p$-order of strong
convergence ($p=1/2$ for Euler scheme) \cite{Milstein1995}:
\begin{gather}
\sup_{t\in [0,\tau_1^{\Lambda})}
\mathbb{E}(|\bfX_{t}-\hat{\bfX}_{t}|^2+|\bfs_{t}-\hat{\bfs}_{t}|^2
~|\mathcal{F}_0)\leq K(1+|\bfs_{0}|^2+|\bfX_{0}|^2)h^{2p},
\end{gather}
and using this inequality together  with (\ref{Appendix::Granwall}),
(\ref{Appendix::jump}) one can come to the conclusion  that the
overall strong error of the scheme is
\begin{gather}
\epsilon_X([0,T))\leq C(\{\tau_n^{\Lambda}\},T)h^{2p},\\
\epsilon_\sigma([0,T))\leq C(\{\tau_n^{\Lambda}\},T)h^{2p}
\end{gather}
where  constant $C$ depends on the sequence of jump times $\{\tau_{n}^{\Lambda}\}$ and length of
the interval $T$. Hence after averaging over all possible jump-times overall error has a strong
asymptotic of the numerical scheme used to propagate the diffusion modes, i.e.
$\overline{\epsilon_X([0,T))}, \overline{\epsilon_\sigma([0,T))} \propto h^{2p}$.
\bibliography{hybrid-scheme.bbl}

\pagebreak\newpage\pagebreak
\section{Figures}


\begin{figure}[!ht]
\begin{center}\vskip 2cm
\includegraphics[angle=-90,scale=0.50]{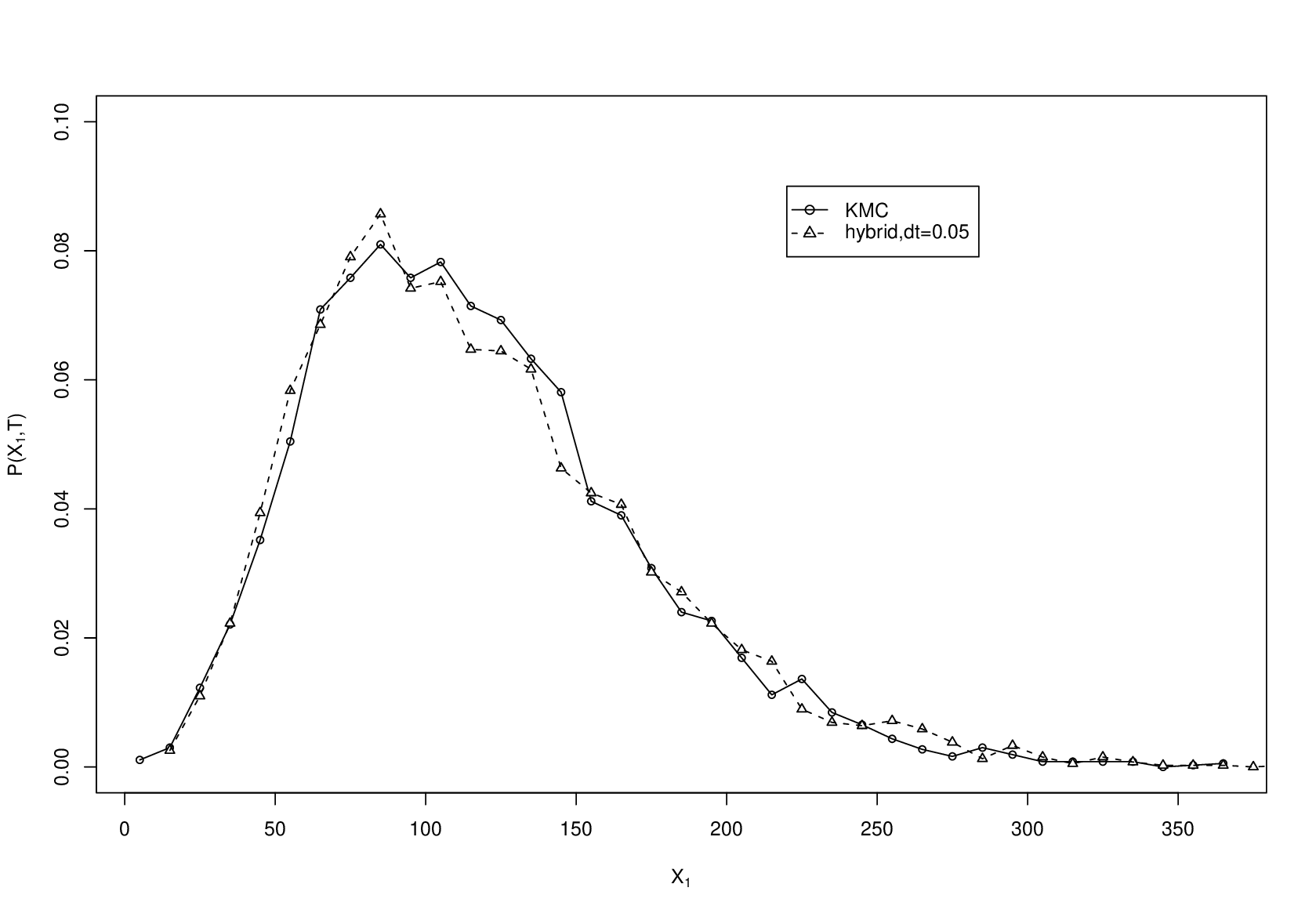}
\end{center}
\caption{\label{Fig::histogramX1} Comparison of empirical stationary distributions $P(X_1,T)$
obtained for $N=10^4$ points for KMC and hybrid scheme with $h=0.1$ at $T=2\times 10^3$ }
\end{figure}
\newpage\pagebreak\newpage
\begin{figure}[!ht]
\begin{center}\vskip 2cm
\includegraphics[angle=-90,scale=0.50]{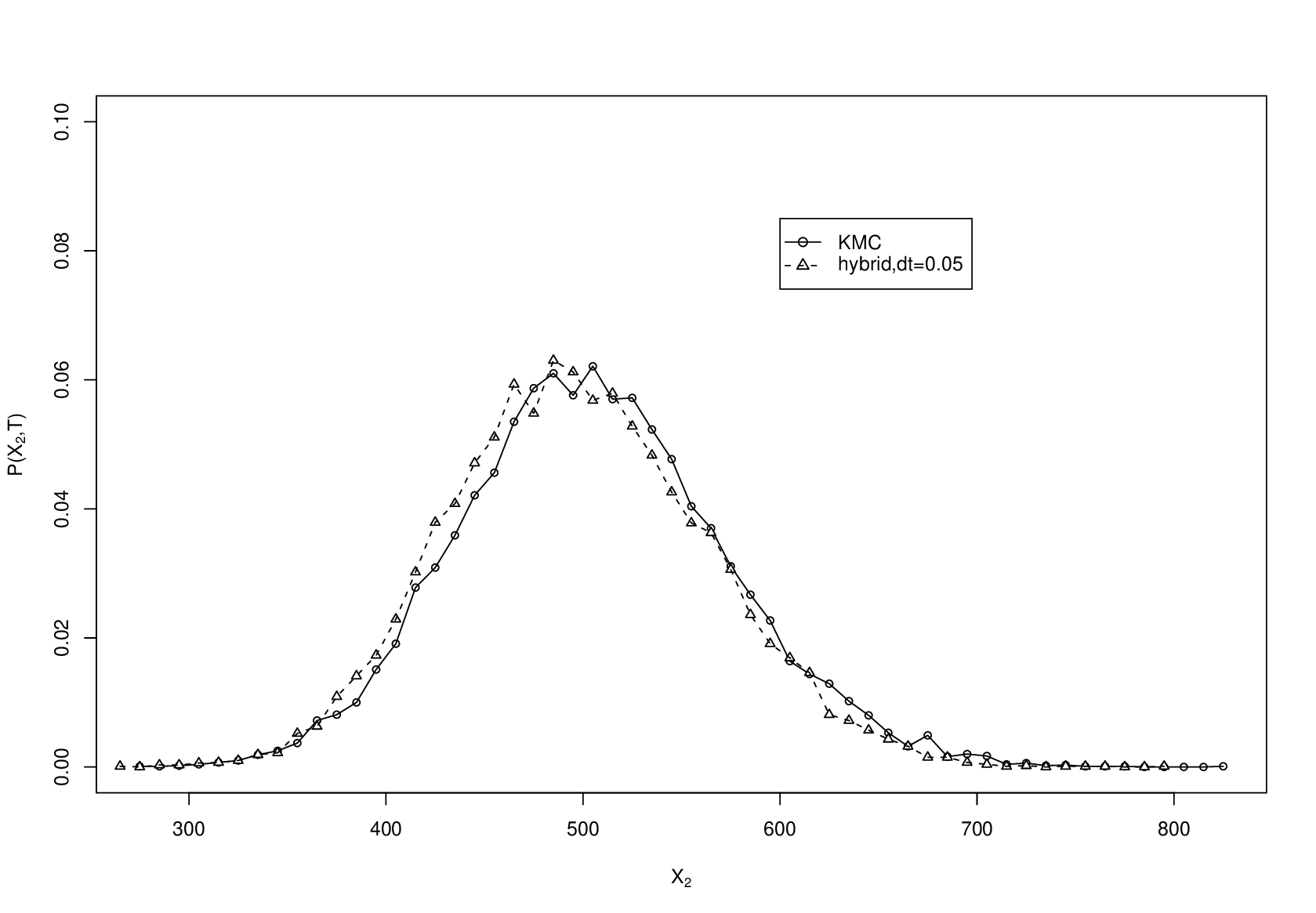}
\end{center}
\caption{\label{Fig::histogramX2} Comparison of empirical stationary distributions $P(X_2,T)$
obtained for $N=10^4$ points for KMC and hybrid scheme with $h=0.1$ at $T=2\times 10^3$ }
\end{figure}
\newpage\pagebreak\newpage
\begin{figure}[!ht]
\begin{center}\vskip 2cm
\includegraphics[angle=-90,scale=0.50]{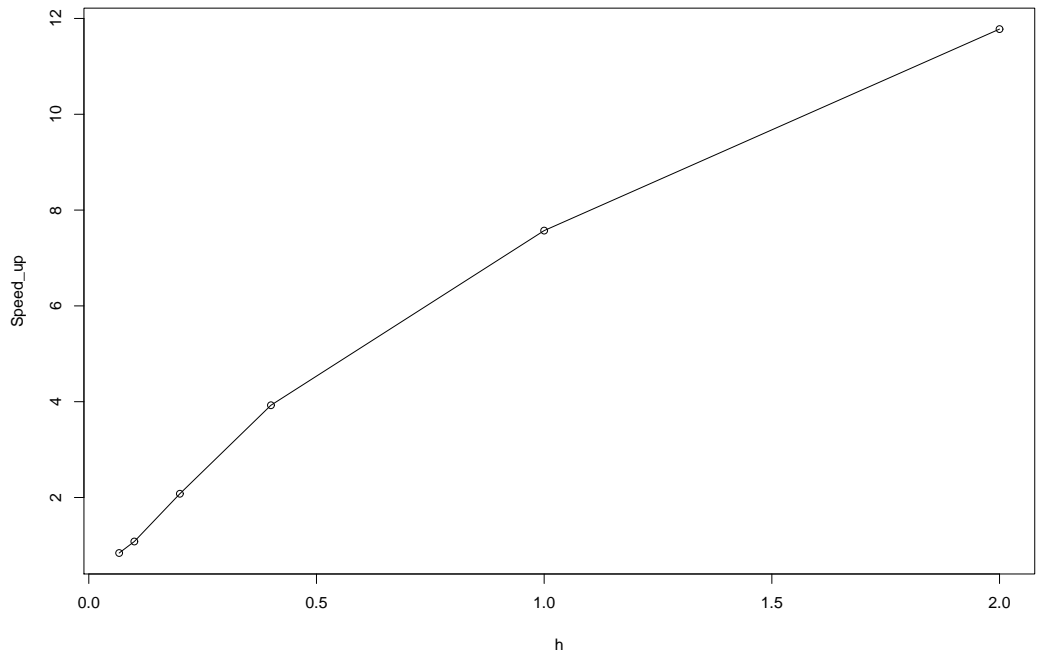}
\end{center}
\caption{\label{Fig::speed_up} Speedup of the hybrid method compared to the  KMC method for
different step-sizes of the grid. }
\end{figure}
}
\end{document}